\documentstyle[11pt,leqno,amsfonts,amssymb]{article}
\newtheorem{definition}{Definition}[section]
\newtheorem{lemma}[definition]{Lemma}
\newtheorem{theorem}[definition]{Theorem}
\newtheorem{proposition}[definition]{Proposition}
\newtheorem{corollary}[definition]{Corollary}
\newtheorem{remark}[definition]{Remark}

\newtheorem{example}[definition]{Example}

 1

\newenvironment{proof}{\noindent{\bf Proof~:}}{\QED\medskip}
\def\QED{\hskip0.1em\hfill\null\ \null\nobreak\hfill
\kern3pt\lower1.8pt\vbox{\hrule\hbox
{\vrule\kern1pt\vbox{\kern1.7pt \hbox{$\scriptstyle
QED$}\kern0.2pt}\kern1pt\vrule}\hrule}}

\newcommand{\ZZ}{{\mathbb{Z}}}
\newcommand{\CC}{{\mathbb{C}}}
\newcommand{\QQ}{{\mathbb{Q}}}

\newcommand{\RR}{{\mathbb{R}}}
\newcommand{\TT}{{\mathbb{T}}}

\newcommand{\CP}{{\mathbb{CP}}}
\newcommand{\wcpn}{{\widetilde{CP}{}^m}}
\newcommand{\slc}{{\frak s}{\frak l}(2,\CC)}

\newcommand{\coker}{\hbox{coker}\,}

\renewcommand{\Im}{\hbox{Im}\,}

\newcommand{\rank}{\hbox{rank}\,}

\newcommand{\inc}{\hookrightarrow}

\newcommand{\PD}{\hbox{PD}}
\newcommand{\la}{\langle}
\newcommand{\ra}{\rangle}

\newcommand{\Lo}{L_{[\omega]}}
\newcommand{\LO}{L_{[\Omega]}}

\begin{document}

\title{\bf Weakly Lefschetz symplectic manifolds}

\author{Marisa Fern\'andez, Vicente Mu\~noz and Luis Ugarte}

\date{February 8, 2005}

\maketitle

\begin{abstract}
The harmonic cohomology of a Donaldson symplectic submanifold and
of an Auroux symplectic submanifold are compared with that of its
ambient space. We also study symplectic manifolds satisfying a
weakly Lefschetz property, that is, the $s$--Lefschetz propery. In
particular, we consider the symplectic blow-ups ${\wcpn}$ of the
complex projective space ${CP}^m$ along weakly Lefschetz
symplectic submanifolds $M\subset {CP}^m$. As an application we
construct, for each even integer $s\geq 2$, compact symplectic
manifolds which are $s$--Lefschetz but not $(s+1)$--Lefschetz.
\end{abstract}

\section{Introduction} \label{introduction}

One of the main results of Hodge theory states that any de Rham
cohomology class on a compact oriented Riemannian manifold has an
unique harmonic representative. In the symplectic setting a notion
of harmonicity can be introduced as follows~\cite{Bry}. Let
$(M,\omega)$ be $2n$--dimensional symplectic manifold. A closed
form $\alpha$ on $M$ is called {\em symplectically harmonic} if
$\delta\alpha=0$, where $\delta$ denotes the Koszul
differential~\cite{Kos}. However, a symplectic version of the
above result does not hold in general. In fact, Mathieu~\cite{Mat}
proved that any de Rham cohomology class has a (not necessarily
unique) symplectically harmonic representative if and only if
$(M,\omega)$ satisfies the hard Lefschetz property, i.e. the map
 \begin{equation} \label{eqn:v1}
  L^{n-k} \colon H^k(M) \longrightarrow H^{2n-k}(M)
 \end{equation}
given by $L^{n-k}[\alpha]=[\alpha\wedge\omega^{n-k}]$ is onto for
all $k \leq n-1$.

In this paper we deal with symplectic manifolds satisfying a
weaker property: following~\cite{FM}, we shall say that
$(M,\omega)$ is an {\it $s$--Lefschetz symplectic manifold}, $0
\leq s\leq n-1$, if (\ref{eqn:v1}) is an epimorphism for all $k
\leq s$. As an obvious fact, whenever $(M,\omega)$ is not hard
Lefschetz, there is some $s\geq 0$ such that $(M,\omega)$ is
$s$--Lefschetz but not $(s+1)$--Lefschetz. So, it seems
interesting to understand the way this phenomenon occurs on
non-hard Lefschetz symplectic manifolds, in particular if there is
some restriction for the possible values of the level $s$ at which
the Lefschetz property can be lost, how this affects to other
symplectic invariants of the manifold, such as the above mentioned
harmonicity, or if the $s$--Lefschetz property is preserved under
the usual constructions of new symplectic manifolds from old ones,
for instance the symplectic blowing up~\cite{McD1}, the Donaldson
symplectic submanifolds~\cite{Do} and the Auroux symplectic
submanifolds~\cite{Au97}. Our purpose in this paper is to explore
these questions, as we explain next.

Regarding symplectic harmonicity, in Section~\ref{someresults} we
recall some results on the harmonic cohomology of $(M,\omega)$ and
show how the $s$--Lefschetz property is related to the existence
problem of symplectically harmonic representatives of de Rham
classes of $M$. Let us denote by $H^{k}_{\rm hr}(M,\omega)$ the
space of {\em harmonic cohomology} in degree $k$, that is, the
subspace of the de Rham cohomology group $H^k(M)$ consisting of
all classes which contain at least one symplectically harmonic
$k$--form. In Proposition~\ref{grados+2} we prove that a
$2n$--dimensional symplectic manifold $(M,\omega)$ is
$s$-Lefschetz if and only if $H^{2n-k}_{\rm hr}(M,\omega) =
H^{2n-k}(M)$ for every $k\leq s$; moreover, the latter condition
implies that $H^{k}_{\rm hr}(M,\omega) = H^k(M)$ for every $k\leq
s+2$. In the proof of this proposition, which can be seen as a
refinement of the result of Mathieu, we follow the approach by
Yan~\cite{Ya} which uses the theory of infinite dimensional
$\slc$--representations.

Section~\ref{Donaldsonsect} is devoted to the study of harmonic
cohomology of Donaldson and Auroux symplectic submanifolds. Given
a compact symplectic manifold $(M,\omega)$ of dimension~$2n$ such
that $[\omega]\in H^2(M)$ admits a lift to an integral cohomology
class, Donaldson proves in~\cite{Do} the existence of a symplectic
submanifold $(Z,\omega_Z)$ of codimension 2 in $M$ which realizes
the Poincar\'e dual of $k\,[\omega]$ for any sufficiently large
integer $k$, and such that the inclusion $\jmath\colon Z
\hookrightarrow M$  is $(n-1)$--connected. We show the following relation between
the harmonic cohomologies $H^*_{\rm hr}(Z,\omega_Z)$ and $H^*_{\rm
hr}(M,\omega)$.

\begin{theorem} \label{harmonicDonaldson}
The inclusion $\jmath\colon Z \hookrightarrow M$ induces an
isomorphism $\jmath^*\colon H^i_{\rm hr}(M,\omega) \longrightarrow
H^i_{\rm hr}(Z,\omega_Z)$ for any $i< n-1$, and a monomorphism for
$i=n-1$. Moreover, $H^i_{\rm hr}(Z,\omega_Z)$ and $H^{i+2}_{\rm
hr}(M,\omega)$ are isomorphic for every $n\leq i\leq 2(n-1)$.
\end{theorem}

Roughly speaking, this result says that a Donaldson symplectic
submanifold inherits essentially the same harmonic cohomology as
that of its ambient space, with the only possible exception of
having less symplectically harmonic forms in the middle degree
$n-1$. Auroux has generalized in~\cite{Au97} Donaldson's
construction. We show that a result like
Theorem~\ref{harmonicDonaldson} does not hold in general for the
Auroux submanifolds. Moreover, the harmonic cohomology of Auroux
symplectic submanifolds has a very different behaviour with
respect to its ambient space, and surprisingly there exist
submanifolds having strictly more harmonic cohomology classes than
their ambient spaces. More concretely, in Example~\ref{example1}
we construct an Auroux submanifold $(Z,\omega_{Z})$ of codimension
2 in a 10-dimensional compact symplectic manifold $(M,\omega)$
such that the inclusion $\jmath\colon Z \hookrightarrow M$ induces
an isomorphism between the de Rham cohomology groups $H^3(Z)$ and
$H^3(M)$, but $\dim H^3_{\rm hr}(Z,\omega_{Z}) >\dim H^3_{\rm
hr}(M,\omega)$.

\medskip

Given a compact symplectic manifold $(M,\omega)$ of dimension
$2n$, we can assume, without loss of generality,  that the
symplectic form $\omega$ is integral (by perturbing and
rescaling). A theorem of Gromov and Tischler~\cite{Gr1, Gr2,Ti}
states that there is a symplectic embedding $i\colon
(M,\omega)\longrightarrow (CP^m,\omega_0)$, with $m\geq {2n+1}$,
where $\omega_0$ is the standard K\"ahler form on ${CP}^m$ defined
by its natural complex structure and the Fubini--Study metric. We
consider the symplectic blow-up ${\wcpn}$ of ${CP}^m$ along the
embedding $i$ (see~\cite{McD1}). Then, ${\wcpn}$ is a simply
connected compact symplectic manifold. In
Section~\ref{symplecticblowups} we study the $s$--Lefschetz
property of ${\wcpn}$, $m\geq 2n+1$. More concretely we have the
following result.

\begin{theorem} \label{thm:main2}
If $(M,\omega)$ is an $s$--Lefschetz compact symplectic manifold
of dimension $2n$, then the symplectic blow-up ${\wcpn}$ $(m\geq
{2n+1})$ is $(s+2)$--Lefschetz. Moreover, if $M$ is parallelizable
and not $s$--Lefschetz then ${\wcpn}$ is not $(s+2)$--Lefschetz.
\end{theorem}

This will be proved in Theorem~\ref{lefschetz-blowup} and
Proposition~\ref{nolefschetz-blowup}. Recently
Cavalcanti~\cite{Cav} has investigated the hard Lefschetz property
of symplectic blow-ups of non-hard Lefschetz symplectic manifolds
along hard Lefschetz symplectic submanifolds; in particular, he
obtains that the symplectic blow-up of a hard Lefschetz symplectic
manifold along a hard Lefschetz submanifold is always hard
Lefschetz. Such a result can be also proved with the arguments of
Theorem~\ref{thm:main2} as we notice in Remark 3.3.

In~\cite{FM} examples of compact symplectic manifolds which are
$s$--Lefschetz but not $(s+1)$--Lefschetz are constructed for each
$s\leq 2$. As an application of Theorem~\ref{thm:main2} and of the
results of Section~\ref{someresults} on the harmonic cohomology of
iterated Donaldson submanifolds of symplectic blow-ups, we prove
in Section~\ref{Donaldsonhr} that for each {\em even\/} integer
number $s\geq 2$, there is a simply connected compact symplectic
manifold of dimension $2(s+2)$ which is $s$--Lefschetz but not
$(s+1)$--Lefschetz. Notice that $2(s+2)$ is the lowest possible
dimension where such a manifold can live. With the same
techniques, we also show a simply-connected symplectic
$10$--manifold which is $3$--Lefschetz but not $4$--Lefschetz.

\section{Harmonic cohomology of $s$-Lefschetz symplectic manifolds} \label{someresults}

We recall some definitions and results about the symplectic
codifferential and symplectically harmonic forms. Let $(M,\omega)$
be a symplectic manifold, that is, $M$ is a differentiable
manifold of dimension~$2n$ and $\omega$ a closed non-degenerate
$2$--form on $M$, the {\it symplectic form}. Denote by
$\Omega^*(M)$, ${\mathfrak X}(M)$ and ${\cal F}(M)$ the algebras
of differential forms, vector fields and differentiable functions
on $M$, respectively. The isomorphism
 $$
 \natural:{\mathfrak X}(M) \longrightarrow \Omega^1(M)
 $$
given by $\natural(X)=\iota_X(\omega)$ for $X\in {\mathfrak
X}(M)$, where $\iota_X$ denotes the contraction by $X$, extends to
an isomorphism of algebras $\natural:\bigoplus_{k\geq 0}{\mathfrak
X}^k(M) \longrightarrow \bigoplus_{k\geq 0}\Omega^k(M)$. Then,
$G=-\natural^{-1}(\omega)$ is the skew-symmetric bivector field
dual to $\omega$. ($G$ is the unique non-degenerate Poisson
structure \cite{Lich} associated with $\omega$.) The Koszul
differential $\delta\colon \Omega^k(M) \longrightarrow
\Omega^{k-1}(M)$ is defined by
 $$
 \delta = [\iota_G,d].
 $$

In \cite{Bry} Brylinski proved that the Koszul differential is a
{\em symplectic codifferential} of the exterior differential with
respect to the {\em symplectic star operator} defined as follows.
Denote  by $\Lambda^k(G)$, $k \geq 0$, the associated pairing
$\Lambda^k(G): {\Omega^k(M)\times\Omega^k(M)}\longrightarrow {\cal
F}(M)$ which is $(-1)^k$--symmetric (i.e.\, symmetric for even
$k$, anti-symmetric for odd $k$). Let $v_M$ be the volume form on
$M$ given by $v_M=\frac {\omega^{n}} {n!}$. Imitating the Hodge
star operator for Riemannian manifolds, the {\em symplectic star
operator}
 $$
 *:\Omega^k(M) \longrightarrow \Omega^{2n-k}(M)
 $$
is defined by the condition $\beta\wedge
(*\alpha)=\Lambda^k(G)(\beta, \alpha)\, v_M$, for $\alpha, \beta
\in \Omega^k(M)$. An easy consequence is that $*^2 = Id$, and if
$\alpha \in \Omega^k(M)$ then
 $$
 \delta(\alpha) =(-1)^{k+1}(*\circ d \circ *)(\alpha).
 $$

Since $\omega$ is a closed form, for any $p,k \geq 0$ the
homomorphism
 $$
 L^p:\Omega^k(M) \longrightarrow \Omega^{2p+k}(M)
 $$
given by $L^p(\alpha)=\alpha\wedge\omega^p$ for $\alpha \in
\Omega^k(M)$, satisfies that $[L^p,d] =L^p \circ d - d \circ L^p =
0$, and therefore it induces a map $L^p: H^k(M) \longrightarrow
H^{2p+k}(M)$ on de Rham cohomology. Relations between the
operators $\iota_G$, $L$, $d$ and $\delta$ are proved by Yan in
\cite{Ya}. Here we shall need the following
 \begin{equation} \label{eqn:v0}
  [L,\delta]=d.
 \end{equation}

A $k$--form $\alpha \in \Omega^k(M)$ is said to be {\em
symplectically harmonic} if $d \alpha=\delta\alpha=0$. Let
$\Omega^k_{\rm hr}(M,\omega)=\{\alpha\in\Omega^k(M) \mid
d\alpha=\delta\alpha=0\}$ be the space of the symplectically
harmonic $k$--forms. Yan proved that for any $k\geq 0$ the map
$L^{n-k}:\Omega^k(M) \longrightarrow \Omega^{2n-k}(M)$ is an
isomorphism. This also induces an isomorphism when restricted to
the subspaces of harmonic forms, as follows from (\ref{eqn:v0}).

\begin{lemma}\hskip-.1cm{\rm \cite{Ya}} (Duality on harmonic forms). \label{dualityforms}
The map
 $$
 L^{n-k} \colon \Omega^{k}_{\rm hr}(M,\omega)\longrightarrow
 \Omega^{2n-k}_{\rm hr}(M,\omega)
 $$
is an isomorphism for $k\geq 0$.
\end{lemma}

For the de Rham cohomology classes of $M$, we consider the vector
space
 $$
 H^k_{\rm hr}(M,\omega)={\Omega^k_{\rm hr}(M,\omega)\over
 \Omega^k_{\rm hr}(M,\omega)\cap {\rm Im}\, d}
 $$
consisting of the cohomology classes in $H^k(M)$ containing at
least one symplectically harmonic form. Lemma~\ref{dualityforms}
implies that the homomorphism
 $$
 L^{n-k}\colon H^{k}_{\rm hr}(M,\omega)\longrightarrow H^{2n-k}_{\rm
 hr}(M,\omega)
 $$
is surjective. (Notice that the duality on harmonic forms may not
be satisfied at the level of the spaces $H^{*}_{\rm
hr}(M,\omega)$.) Since $H^{2n-k}_{\rm hr}(M,\omega)$ is a subspace
of the de Rham cohomology $H^{2n-k}(M)$, we conclude (see
\cite[Corollary 1.7]{IRTU})
\begin{equation}\label{highdegree}
 H^{2n-k}_{\rm hr}(M,\omega) = \Im(L^{n-k} \colon H^{k}_{\rm
 hr}(M,\omega)\longrightarrow H^{2n-k}(M)).
\end{equation}

A nonzero $k$--form $\alpha$, with $k \leq n$,  is called {\em
effective} if $L^{n-k+1}(\alpha)=0$. A cohomology class $a\in
H^k(M)$ is said to be {\em primitive\/} if
  $L^{n-k+1}(a)=0$ in $H^{2n-k+2}(M)$.
In~\cite[page 46]{LM} the following result is proved.

\begin{lemma}\hskip-.1cm{\rm \cite{LM}}. \label{effectiveforms}
If $\alpha$ is an effective $k$--form, then there is a constant
$c$ such that its symplectic star operator  $*\alpha$ satisfies
$*\alpha = c\, L^{n-k}(\alpha)$.
\end{lemma}

Therefore, any closed effective $k$--form $\alpha$ on $(M,\omega)$
is symplectically harmonic because $d *\alpha = c\,\,d\,
L^{n-k}(\alpha)=0$ by Lemma~\ref{effectiveforms}. In particular,
every closed $1$--form is symplectically harmonic since it is
effective, so $H^{1}_{\rm hr}(M,\omega)=H^1(M)$. Moreover, we have

\begin{proposition} \label{gradop}
Let $(M,\omega)$ be a symplectic manifold of dimension $2n$.
Suppose that there exists some integer $k\leq n$ with
$H^{2n-k+2}(M)=0$. Then, for  any closed $k$--form $\alpha \in
\Omega^k(M)$, there is a closed $k$--form $\widetilde{\alpha}$
such that $\widetilde{\alpha}$ is cohomologous to $\alpha$ and is
symplectically harmonic. In particular, $H^{k}_{\rm
hr}(M,\omega)=H^k(M)$.
\end{proposition}

\begin{proof}
Let $a=[\alpha] \in H^k(M)$. We will find a symplectically
harmonic representative of the cohomology class $a$. Since
$L^{n-k+1}(\alpha)$ is a closed $(2n-k+2)$--form and
$H^{2n-k+2}(M)$ is zero, there is some $\beta\in
\Omega^{2n-k+1}(M)$ such that $L^{n-k+1}(\alpha)=d\beta$. But the
map $L^{n-k+1}\colon \Omega^{k-1}(M)\longrightarrow
\Omega^{2n-k+1}(M)$ is surjective, so there exists
$\gamma\in\Omega^{k-1}(M)$ satisfying $\beta=L^{n-k+1}(\gamma)$.
Hence $L^{n-k+1}(\alpha)=d\beta=L^{n-k+1}(d\gamma)$, i.e.,\
$L^{n-k+1}(\alpha-d\gamma)=0$.

Consider $\widetilde{\alpha}=\alpha-d\gamma$, which is
cohomologous to $\alpha$. Using Lemma~\ref{effectiveforms},
$L^{n-k+1}(\widetilde{\alpha})=0$ implies that $*
\widetilde{\alpha}=c\, L^{n-k}(\widetilde{\alpha})$ for some
constant $c$. Thus $d* \widetilde{\alpha}=c\,
L^{n-k}(d\widetilde{\alpha})=0$ and the $k$--form
$\widetilde{\alpha}$ is symplectically harmonic.
\end{proof}

For the de Rham classes in $H^2(M)$,  Mathieu proved the following
result.

\begin{lemma}\hskip-.1cm{\rm \cite{Mat}}. \label{2forms}
Any cohomology class of degree $2$ has a symplectically harmonic
representative.
\end{lemma}

As a consequence of the previous results, if $(M,\omega)$ is a
simply connected  compact symplectic manifold then every class in
$H^k(M)$ has a symplectically harmonic representative for~$k \leq
3$.

\begin{corollary} \label{6simplyconnected}
 Let $(M,\omega)$ be a compact simply connected symplectic
 manifold of dimension~$6$. Then every de Rham cohomology class of degree
 $k \neq 4$ admits a symplectically harmonic representative.
\end{corollary}

\medskip
Recall that a symplectic manifold $(M,\omega)$ of
dimension $2n$ is said to be {\em $s$--Lefschetz\/} with $0 \leq
s\leq n-1$, if the map $L^{n-k} \colon H^k(M) \longrightarrow
H^{2n-k}(M)$ is an epimorphism for all $k \leq s$. In the compact
case we actually have that $L^{n-k}$ are isomorphisms because of
Poincar\'e duality. Note that $M$ is $(n-1)$--Lefschetz if $M$
satisfies the hard Lefschetz theorem.

\begin{proposition} \label{grados+2}
 Let $(M,\omega)$ be a symplectic manifold of dimension
 $2n$ and let $s\leq n-1$. Then the following statements are
 equivalent:
 \begin{enumerate}
 \item $(M,\omega)$ is $s$--Lefschetz.
\item $H^{k}_{\rm hr}(M,\omega) = H^k(M)$ for every
 $k\leq s+2$, and $H^{2n-k}_{\rm hr}(M,\omega) = H^{2n-k}(M)$ for every $k\leq s$.
\item $H^{2n-k}_{\rm hr}(M,\omega) = H^{2n-k}(M)$ for every $k\leq
s$.
 \end{enumerate}
\end{proposition}

\begin{proof}
Clearly $(ii)$ implies $(iii)$. Let us see also that $(iii)$
implies $(i)$. Let $k\leq s$. From~(\ref{highdegree}) we have that
 $$
 H^{2n-k}_{\rm hr}(M,\omega) = \Im(L^{n-k}\mid_{H^{k}_{\rm
 hr}(M,\omega)} \colon H^{k}_{\rm
 hr}(M,\omega) \inc H^k(M) \longrightarrow H^{2n-k}(M)).
 $$
If $H^{2n-k}_{\rm hr}(M,\omega)= H^{2n-k}(M)$, then the map
$L^{n-k}\mid_{H^{k}_{\rm hr}(M,\omega)}$ is onto, and therefore
the homomorphism $L^{n-k}\colon H^{k}(M)\longrightarrow
H^{2n-k}(M)$ must also be onto. So $M$ is $s$--Lefschetz.

We want to show that $(i)$ implies $(ii)$. It is enough to prove
that $H^{k}_{\rm hr}(M,\omega) = H^k(M)$ for every $k\leq s+2$,
because in this case, for $k\leq s$, we have $H^{2n-k}_{\rm
hr}(M,\omega) = \Im(L^{n-k} \colon H^{k}(M)\longrightarrow
H^{2n-k}(M))= H^{2n-k}(M)$ using the $s$--Lefschetz property.

Let us see that $H^{k}_{\rm hr}(M,\omega) = H^k(M)$ for every
$k\leq s+2$, by induction on $s$. For $s=0$, we recall that $M$ is
$0$--Lefschetz as this is satisfied by every symplectic manifold.
Now for any symplectic manifold, any class of degree $1$ admits an
harmonic representative by Lemma~\ref{effectiveforms}, and we also
know that any class of degree $2$ admits an harmonic
representative by Lemma~\ref{2forms}.

Now take $s>0$, and suppose that if $(M,\omega)$ is
$(s-1)$--Lefschetz, it holds $H^{k}_{\rm hr}(M,\omega) = H^k(M)$
for $k\leq s+1$. We have to prove that $H^{s+2}_{\rm hr}(M,\omega)
= H^{s+2}(M)$ if $M$ is $s$--Lefschetz. Let $\alpha$ be a closed
element of degree $s+2$. Consider the map
$L^{n-s-1}:\Omega^{s+2}(M) \longrightarrow \Omega^{2n-s}(M)$. Then
$L^{n-s-1}(\alpha)$ is a closed $(2n-s)$--form. By the
$s$--Lefschetz property there is a closed $s$--form $h$ (which we
may suppose to be symplectically harmonic, by induction
hypothesis) such that
 $$
 L^{n-s-1}(\alpha) = L^{n-s} (h) + d\beta\, ,
 $$
for some $\beta \in \Omega^{2n-s-1}(M)$. By the surjectivity of
$L^{n-s-1}\colon \Omega^{s+1}(M)\longrightarrow
\Omega^{2n-s-1}(M)$ we get the existence of some $(s+1)$--form
$\gamma$ with $\beta= L^{n-s-1}(\gamma)$. Therefore
$L^{n-s-1}(\alpha) = L^{n-s}(h) + L^{n-s-1} (d\gamma)$ and hence
 \begin{equation} \label{e}
 L^{n-s-1}( \alpha - L(h) - d\gamma) =0.
 \end{equation}

Put $\widetilde{\alpha}=\alpha - L(h) - d\gamma$. By (\ref{e}) and
Lemma  \ref{effectiveforms} we have that $* \,\widetilde{\alpha}=
c\, L^{n-s-2}(\widetilde{\alpha})$ for some constant $c$.
Therefore $\widetilde{\alpha}$ is symplectically harmonic. On the
other hand, since $h$ is symplectically harmonic we see that
$L(h)$ is symplectically harmonic using (\ref{eqn:v0}). Hence
$\alpha - d\gamma$ is symplectically harmonic and cohomologous to
the original $\alpha$.
\end{proof}

Notice that this result implies that every de Rham cohomology
class of $M$ admits a symplectically harmonic representative if
and only if $(M,\omega)$ is hard Lefschetz, which is Mathieu's
theorem. Also, if $M$ is a simply connected compact symplectic
manifold of dimension~6 then $(M,\omega)$ is hard Lefschetz if and
only if every cohomology class of degree $4$ has a symplectically
harmonic representative.

\bigskip

If $M$ is a manifold of finite type, i.e. all the de Rham
cohomology groups $H^k(M)$ are finite dimensional, then we shall
denote by $b_{k}^{\rm hr}(M,\omega)$ the dimension of the space
$H^{k}_{\rm hr}(M,\omega)$. As usual, the Betti numbers of $M$
will be denoted by $b_k(M)=\dim H^k(M)$.

It is well-known that if $(M,\omega)$ is compact and hard
Lefschetz, the odd Betti numbers of $M$ are even. When
$(M,\omega)$ is $s$--Lefschetz we have the following proposition.

\begin{proposition}\label{Bettinumbers}
Let $(M,\omega)$ be a compact symplectic manifold of dimension
$2n$. Suppose that $(M,\omega)$ is $s$--Lefschetz with $s\leq
{n-1}$. Then  the odd Betti numbers $b_{2i-1}(M)$ are even for
${2i-1} \leq s$, and $b_{2n-2j+1}^{\rm hr}(M,\omega)$ is even for
$s < 2j-1 \leq s+2$.
\end{proposition}

\begin{proof}
Put $k=2i-1\leq s$.  Let us consider the non-singular pairing
 $$
 p\colon H^k(M)\otimes H^{2n-k}(M) \longrightarrow \ {\RR}
 $$
given by
 $$
 p([\alpha],[\beta])=\int_M \alpha\wedge\beta,
 $$
for $[\alpha]\in H^k(M)$ and $[\beta]\in H^{2n-k}(M)$. Let
$\langle\ ,\ \rangle$ be the skew-symmetric bilinear form defined
on  $H^k(M)$ by
 $$
 \langle [\alpha] ,[\alpha'] \rangle= p([\alpha], L^{n-k}[\alpha']),
 $$
for $[\alpha] ,[\alpha'] \in H^k(M)$. The rank of  $\langle\ ,\
\rangle$ is an even number (see page 4, \cite{LM}). The
non-singularity of $p$ implies that the rank of $\langle\ ,\
\rangle$ equals the rank of the map $L^{n-k} \colon
H^k(M)\longrightarrow H^{2n-k}(M)$, that is, $\rank{\langle\ ,\
\rangle}=b_{2n-k}(M)$ since $(M,\omega)$ is $s$--Lefschetz. Hence
$b_k(M)$ is even by Poincar\'e duality.

For the final part, take $k=2j-1$ with $s<k \leq {s+2}$. Now, the
previous argument also shows that $b_{2n-k}^{\rm hr}(M,\omega)$ is
even because the $s$--Lefschetz property implies $H^k(M)=
H^{k}_{\rm hr}(M,\omega)$ by Proposition~\ref{grados+2} and, on
other hand, $H^{2n-k}_{\rm hr}(M,\omega)=\Im(L^{n-k})$, therefore
the rank of  $\langle\ ,\ \rangle$ is an even number which equals
$\dim \Im(L^{n-k})=b_{2n-k}^{\rm hr}(M,\omega)$.
\end{proof}

\section{Harmonic cohomology of Donaldson and Auroux symplectic submanifolds}
\label{Donaldsonsect}

In this section we study the relation between the harmonic
cohomology of Donaldson and Auroux symplectic submanifolds and that of the
ambient space.

Recall that given a compact symplectic manifold $(M,\omega)$ of
dimension $2n$ such that $[\omega]\in H^2(M)$ admits a lift to an
integral cohomology class, Donaldson proves~\cite{Do} the
existence of a symplectic submanifold $Z$ of codimension 2 in $M$
that realizes the Poincar\'e dual of $k\, [\omega]$ for any
sufficiently large integer $k$. Moreover, the inclusion
$\jmath\colon Z \hookrightarrow M$ is $(n-1)$--connected, that is,
$\jmath^*\colon H^i(M) \longrightarrow H^i(Z)$ is an isomorphism
for $i<(n-1)$, and a monomorphism for $i=(n-1)$. Let us denote by
$\omega_Z=\jmath^*\omega$ the symplectic form on $Z$.

\bigskip
\noindent{\bf Proof of Theorem~\ref{harmonicDonaldson}~:} We use
here the following property given in~\cite[Lemma 4.3]{Yam}: If
$(M,\omega)$ is a symplectic manifold of dimension~$2n$, then for
any $2\leq i\leq n$ the subspace $H^{i}_{\rm hr}(M,\omega)$ of
$H^{i}(M)$ is given by
 \begin{equation}\label{YamadaM}
 H^{i}_{\rm hr}(M,\omega) = P_{i}(M,\omega) + L_{[\omega]}\left(
 H^{i-2}_{\rm hr}(M,\omega) \right) ,
 \end{equation}
where $P_{i}(M,\omega)=\{\, a\in H^{i}(M) \ \mid \
L_{[\omega]}^{n-i+1} (a)=0 \}$.

Similarly, for the Donaldson symplectic submanifold $Z$ we have
$$
H^{i}_{\rm hr}(Z,\omega_Z) = P_{i}(Z,\omega_Z) + L_{[\omega_Z]}
\left( H^{i-2}_{\rm hr}(Z,\omega_Z) \right) ,
$$
for any $2\leq i\leq n-1$, where $P_{i}(Z,\omega_Z)=\{\, b\in
H^{i}(Z) \ \mid \ L_{[\omega_Z]}^{n-i} (b)=0 \}$.

On the other hand, in~\cite{FM} it is proved that for any $i\geq
n$, a cohomology class $a\in H^i(M)$ satisfies $\jmath^* a=0$ if
and only if $a\cup [\omega]=0$.

Let us prove first that $\jmath^*(P_{i}(M,\omega)) \subset
P_{i}(Z,\omega_Z)$ for any $2\leq i\leq n-1$. Given $a\in
P_{i}(M,\omega)$, let us consider $b=\jmath^* a \in H^i(Z)$. Since
$0=L_{[\omega]}^{n-i+1} (a) = a\cup [\omega]^{n-i+1}$, and
$n+1\leq 2n-i$ (because $n-1\geq i$), the cohomology class
$L_{[\omega]}^{n-i} a \in H^{2n-i}(M)$ satisfies $\jmath^*
(L_{[\omega]}^{n-i} a) =0$. But $\jmath^*\circ L_{[\omega]} =
L_{[\omega_Z]} \circ \jmath^*$, which implies
$L^{n-i}_{[\omega_Z]}(b)=\jmath^*(L^{n-i}_{[\omega]}a)=0$, that
is, $b\in P_{i}(Z,\omega_Z)$. Now it is easy to see that
$\jmath^*\colon P_{i}(M,\omega) \longrightarrow P_{i}(Z,\omega_Z)$
is an isomorphism for $i<(n-1)$ and a monomorphism for $i=(n-1)$,
because $\jmath^*\colon H^{i}(M) \longrightarrow H^{i}(Z)$ is.

Now, we prove by induction that $\jmath^*(H^i_{\rm hr}(M,\omega))
\subset H^i_{\rm hr}(Z,\omega_Z)$ for any $i\leq (n-1)$. This is
clear for $i=0,1$, because $H^i_{\rm hr} = H^i$. Let us fix $i$
with $2\leq i\leq (n-1)$, and suppose that the inclusion holds in
any degree $<i$. Since $\jmath^*\circ L_{[\omega]} =
L_{[\omega_Z]} \circ \jmath^*$, and $\jmath^*$ takes the primitive
classes of degree $i$ on $M$ to primitive classes of degree $i$ on
the submanifold $Z$, the induction hypothesis and~(\ref{YamadaM})
imply that
$$
\begin{array}{rl}
\jmath^*(H^i_{\rm hr}(M,\omega)) \!\!\! & =  \jmath^*
(P_{i}(M,\omega)) + L_{[\omega_Z]}\left( \jmath^* (H^{i-2}_{\rm
hr}(M,\omega))
\right) \\[9pt]
& \subset  P_{i}(Z,\omega_Z) + L_{[\omega_Z]}\left( H^{i-2}_{\rm
hr}(Z,\omega_Z) \right) \\[8pt]
& = H^i_{\rm hr}(Z,\omega_Z) .
\end{array}
$$
Therefore, for any $i\leq (n-1)$, we have the map $\jmath^*\colon
H^i_{\rm hr}(M,\omega) \longrightarrow H^i_{\rm hr}(Z,\omega_Z)$,
which is just the restriction to the space of harmonic cohomology
classes of the homomorphism $\jmath^* \colon H^i(M)
\longrightarrow H^i(Z)$. Thus, $\jmath^*$ is injective for $i\leq
(n-1)$. Finally, an inductive argument as above allows us to
conclude that $\jmath^*$ is surjective for $i< (n-1)$.

To complete the proof, it remains to see that $H^i_{\rm
hr}(Z,\omega_Z)$ and $H^{i+2}_{\rm hr}(M,\omega)$ are isomorphic
for every $n\leq i\leq 2(n-1)$. Let us consider the spaces $A$ and
$B$ given by
$$
\begin{array}{rl}
A \!\!\! & = \ker ( L^{i-n+2}_{[\omega]}\colon H^{2n-i-2}_{\rm
hr}(M,\omega) \longrightarrow H^{i+2}(M) ),\\[9pt]
B \!\!\! & = \ker (L^{i-n+1}_{[\omega_Z]}\colon H^{2n-i-2}_{\rm
hr}(Z,\omega_Z) \longrightarrow H^{i}(Z) ),
\end{array}
$$
where $n\leq i\leq 2n-2$. Next we see that $\jmath^*$ induces an
isomorphism between $A$ and $B$. Given $a\in A$, we denote
$b=\jmath^*(a) \in H^{2n-i-2}(Z)$. Since $2n-i-2<n-1$, from the
first part of the proof it follows that $b\in H^{2n-i-2}_{\rm
hr}(Z,\omega_Z)$. Moreover, since $a\cup[\omega]^{i-n+2}=0$ if and
only if
$\jmath^*(a\cup[\omega]^{i-n+1})=b\cup[\omega_Z]^{i-n+1}=0$, we
have that $b\in B$, that is, $\jmath^* (A) \subset B$. Again, from
the first part of the proof we conclude that the map $\jmath^*
\colon A \longrightarrow B$ is an isomorphism, because
$(2n-i-2)<(n-1)$.

Finally, as an immediate consequence of~(\ref{highdegree}) we get
$$
\begin{array}{rl}
H^{i}_{\rm hr}(Z,\omega_Z) \!\!\! & = \Im (
L^{i-n+1}_{[\omega_Z]}\colon
H^{2n-i-2}_{\rm hr}(Z,\omega_Z) \longrightarrow H^{i}(Z))
\ \cong\ H^{2n-i-2}_{\rm hr}(Z,\omega_Z) / B \\[9pt]
& \cong \ H^{2n-i-2}_{\rm hr}(M,\omega) / A \ \cong\ \Im (
L^{i-n+2}_{[\omega]}\colon H^{2n-i-2}_{\rm hr}(M,\omega)
\longrightarrow H^{i+2}(M) ) \\[9pt]
& = H^{i+2}_{\rm hr}(M,\omega),
\end{array}
$$
for any $n\leq i\leq (2n-2)$, so $b_i^{\rm
hr}(Z,\omega_Z)=b_{i+2}^{\rm hr}(M,\omega)$ for any such $i$.
\QED\medskip

\medskip

From now on, by an {\it iterated} Donaldson symplectic submanifold
$(Z_l,\omega_l)$ of $(M,\omega)$ we shall mean a symplectic
manifold obtained as
$$
(Z_l,\omega_l)\subset (Z_{l-1},\omega_{l-1}) \subset \cdots
\subset (Z_1,\omega_1)\subset (Z_0=M,\omega_0=\omega),
$$
where $(Z_i,\omega_i)$ is a Donaldson symplectic submanifold of
$(Z_{i-1},\omega_{i-1})$, for any $1\leq i \leq l$.

\begin{corollary}\label{q-submanifold-2}
If $(Z_l,\omega_l)$ is an iterated Donaldson symplectic
submanifold of $(M^{2n},\omega)$, then $b_{n-l}^{\rm
hr}(Z_l,\omega_l) \geq b_{n-l}^{\rm hr}(M,\omega)$ and
$$
\begin{array}{ll}
&  b_i(Z_l) - b_i^{\rm hr}(Z_l,\omega_l)= b_i(M) - b_i^{\rm
hr}(M,\omega)\, , \qquad \mbox{ for } i\leq n-l-1 \, , \\[12pt]
&  b_i(Z_l) - b_i^{\rm hr}(Z_l,\omega_l)= b_{i+2l}(M) -
b_{i+2l}^{\rm hr}(M,\omega)\, , \qquad \mbox{ for } i\geq n-l+1 \,
.
\end{array}
$$
\end{corollary}

\begin{proof}
Applying $l$ times Theorem~\ref{harmonicDonaldson} we have
$b_i^{\rm hr}(Z_l,\omega_l)=b_{i+2l}^{\rm hr}(M,\omega)$ for any
$i\geq n-l+1$. Since $b_{2n-2l-i}(Z_l)=b_{2n-2l-i}(M)$, the
Poincar\'e duality for $Z_l$ and $M$ implies that
$b_{i}(Z_l)=b_{i+2l}(M)$. This proves the corollary for any $i\geq
n-l+1$. For the remaining values of $i$, the result follows
directly from Theorem~\ref{harmonicDonaldson}.
\end{proof}

\medskip
Next we want to show that a result like
Theorem~\ref{harmonicDonaldson} for the Auroux submanifolds does
not hold in general.

Suppose that $(M,\omega)$ is a compact symplectic manifold of
dimension $2n$ with $[\omega]\in H^2(M)$ admitting a lift to an
integral cohomology class, and let $E$ be any hermitian vector
bundle over $M$ of rank $r$. Then, in \cite{Au97} Auroux
constructs symplectic submanifolds $(Z_{r},\omega_{Z_r}) \inc
(M,\omega)$ of dimension $2(n-r)$ whose Poincar\'e duals are
$\PD[Z_r]=c_r(E\otimes L^{\otimes k}) = k^r[\omega]^r +
k^{r-1}c_1(E)[\omega]^{r-1} +\cdots + c_r(E)$ for any integer
number $k$ large enough, where we denote by $c_i(E)$ the $i^{th}$
Chern class of the vector bundle $E$, and by $L$ the complex line
bundle over $M$ with first Chern class $c_1(L)=[\omega]$. These
submanifolds also satisfy a Lefschetz theorem on hyperplane
sections, that is, the inclusion $\jmath\colon Z_{r}\inc M$
induces $\jmath^*\colon H^i(M)\to H^i(Z_{r})$ which is an
isomorphism for $i<(n-r)$ and a monomorphism for $i=(n-r)$.

The strongest result in the direction of
Theorem~\ref{harmonicDonaldson} for the Auroux submanifolds
follows from ~\cite[Theorem 4.4]{FMS}. There it is proved that,
for an Auroux submanifold $Z_{r}\inc M$, for large enough~$k$, and
for each $s\leq (n-r-1)$, if $M$ is $s$--Lefschetz then $Z_{r}$ is
also $s$--Lefschetz. In this situation, we have, thanks to
Proposition~\ref{grados+2}, that $H^i_{\rm hr}(M,\omega) \cong
H^i(M)$ and $H^i_{\rm hr}(Z_r,\omega_{Z_r}) \cong H^i(Z_r)$, for
$i\leq s+2$. Therefore it follows that there is an isomorphism
$\jmath^*\colon H^i_{\rm hr}(M,\omega) \longrightarrow H^i_{\rm
hr}(Z_{r},\omega_{Z_r})$, for any $i\leq \min\{s+2,n-r-1\}$, and a
monomorphism in the case $i=(n-r) \leq (s+2)$.

To disprove a result like Theorem~\ref{harmonicDonaldson} for
Auroux submanifolds, we shall see examples of different behaviours
in the simplest case, i.e.,\ when $M$ is not $1$--Lefschetz. By
the above, $H^i_{\rm hr}(M,\omega) \cong H^i_{\rm
hr}(Z_r,\omega_{Z_r})$, for $i=1,2$. So the first case to look at is the
study of the relation between
 $$
 H^3_{\rm hr}(M,\omega) \quad \hbox{ and } \quad H^3_{\rm hr}(Z_r,\omega_{Z_r}).
 $$
In general, to compare them, we are going to assume $n-r>3$, so
that there is an isomorphism $\jmath^*:H^3(M)\to H^3(Z_r)$. We
need the following lemma.

\begin{lemma} \label{lem:hr3}
 Suppose that $(Z_r,\omega_{Z_r})\inc (M,\omega)$ is an Auroux
 symplectic submanifold, and $n-r>3$.
In the situation above,
 \begin{enumerate}
\item $b_3^{\rm hr}(M,\omega)= b_3(M) +\dim \ker \left(
 \Lo^{n-2}:H^1(M) \to H^{2n-3}(M)\right) \\  -
  \dim \ker \left(
 \Lo^{n-1}:H^1(M) \to H^{2n-1}(M)\right)$,
\item $b_3^{\rm hr}(Z_r,\omega_{Z_r})= b_3(M) +\dim \ker \left(
 \Lo^{n-r-2}\cup c_r(E\otimes L^{\otimes k}):H^1(M) \to H^{2n-3}(M)\right) \\
- \dim \ker \left( \Lo^{n-r-1}\cup c_r(E\otimes L^{\otimes k}):H^1(M) \to
H^{2n-1}(M)\right)$,

 \noindent where $\cup c_r(E\otimes L^{\otimes k}): H^*(M) \to H^{*+2r}(M)$
is interpreted as a map in cohomology.
\end{enumerate}
\end{lemma}

\begin{proof}
Let us start by computing $H^3_{\rm hr}(M,\omega)$. By
(\ref{YamadaM}),
 $$
 H^{3}_{\rm hr}(M,\omega) = P_{3}(M,\omega) + \Lo\left(
 H^{1}_{\rm hr}(M,\omega) \right) ,
 $$
where $P_{3}(M,\omega)=\{\, a\in H^3(M) \ \mid \ \Lo^{n-2} (a)=0
\}$. In the case $i=1$, we have that $H^{1}_{\rm
hr}(M,\omega)=H^1(M)$. Clearly
 $$
 P_{3}(M,\omega) \cap \Lo\left( H^{1}(M) \right) =
 \Lo \left( \ker(\Lo^{n-1}: H^1(M) \to
 H^{2n-1}(M))\right).
 $$
On the other hand, $P_{3}(M,\omega)=\ker \left(\Lo^{n-2}\colon
H^3(M) \to H^{2n-1}(M)\right)$ is dual, via Poincar\'e duality, to
\, $\coker \left(\Lo^{n-2} \colon H^1(M) \to H^{2n-3}(M)\right)$.
Therefore
 \begin{eqnarray*}
  b_3^{\rm hr}(M,\omega) &=& \dim \coker \left(\Lo^{n-2}
\colon H^1(M) \to H^{2n-3}(M)\right) \\ && + \dim \Lo (H^1(M)) -
 \dim \Lo \left( \ker(\Lo^{n-1}: H^1(M) \to
 H^{2n-1}(M))\right) \\
  &=& b_3(M) - b_1(M) + \dim \ker \left(\Lo^{n-2}
 \colon H^1(M) \to H^{2n-3}(M)\right) \\ && + b_1(M)-
 \dim \ker \left(\Lo^{n-1}: H^1(M) \to H^{2n-1}(M)\right) \\
 &=& b_3(M) +\dim \ker \left(
   \Lo^{n-2}:H^1(M) \to H^{2n-3}(M)\right) \\ && - \dim \ker \left(
  \Lo^{n-1}:H^1(M) \to H^{2n-1}(M)\right).
 \end{eqnarray*}
This proves {\it (i)}. Now we move on to compute $H^3_{\rm
hr}(Z_r,\omega_{Z_r})$. First, note that for $i<n-r$, if $a\in
H^{2n-2r-i}(M)$, we have that
 $$
 j^*(a)=0 \iff a\cup c_r(E\otimes L^{\otimes k})=0.
 $$
Certainly, $j^*(a)=0$ is equivalent to
 $$
  0 = \int_{Z_r} j^*(a) \cup j^*(b) =
  \int_M a \cup b \cup c_r(E\otimes L^{\otimes k}),
 $$
for any $b\in H^{i}(M)\cong H^{i}(Z_r)$. We use that $\PD[Z_r]=
c_r(E\otimes L^{\otimes k})$ for the second inequality. This is
equivalent to $a\cup c_r(E\otimes L^{\otimes k})=0$. With the aid
of this, and using {\it (i)}, we have
\begin{eqnarray*}
  b_3^{\rm hr}(Z_r,\omega_{Z_r}) &=& b_3(Z_r) +\dim \ker \left(
   L_{[\omega_{Z_r}]}^{n-r-2}:H^1(Z_r) \to H^{2n-2r-3}(Z_r)\right) \\ && - \dim \ker \left(
   L_{[\omega_{Z_r}]}^{n-r-1}:H^1(Z_r) \to H^{2n-2r-1}(Z_r)\right) \\
   &=& b_3(M) +\dim \ker \left(
  \Lo^{n-r-2}\cup c_r(E\otimes L^{\otimes k}):H^1(M) \to H^{2n-3}(M)\right) \\ && - \dim \ker \left(
   \Lo^{n-r-1}\cup c_r(E\otimes L^{\otimes k}):H^1(M) \to
   H^{2n-1}(M)\right).
 \end{eqnarray*}
\end{proof}

Next we exhibit examples of compact symplectic manifolds
$(X,\Omega)$ having Auroux submanifolds $(Z_r,\Omega_{Z_r})$ such
that $b_3^{\rm hr}(Z_r,\Omega_{Z_r}) \not= b_3^{\rm
hr}(X,\Omega)$. To define $X$, first we consider the simply
connected nilpotent Lie group $G$ of dimension $6$ consisting of
all the matrices of the form
 $$
 \pmatrix{1 & y & t+z & \frac{t}{2} & u + \frac{y^2}{2} & v \cr
 0& 1 & x & \frac{x}{2} & y + \frac{x^2}{2} & xy + \frac{x^3}{6} \cr
 0&0& 1 &0&0& y \cr
 0&0&0& 1 & 2x & x^2 \cr
 0&0&0&0&1& x \cr
 0&0&0&0&0& 1 \cr},
 $$
where $x, y, z, t, u, v \in \RR$. With respect to this global
system of coordinates, the forms
    $$
    \alpha_1=dx,\ \alpha_2=dy,\  \alpha_3=dz,\
    \alpha_4=dt - ydx,\ \alpha_5=du - tdx,\ \alpha_6=dv - (z+t)dy
    - \left(u + \frac{y^2}{2}\right) dx
    $$
constitute a basis of left invariant $1$--forms on $G$, and they
satisfy
 $$
 d\alpha_1 = d\alpha_2 = d\alpha_3 =0, \quad
  d\alpha_4 =
 \alpha_{12}, \quad
 d\alpha_5 = \alpha_{14}, \quad
 d\alpha_6 =
 \alpha_{15} + \alpha_{23} + \alpha_{24},
 $$
where we denote $\alpha_{ij\cdots
k}=\alpha_i\wedge\alpha_j\wedge\cdots \wedge\alpha_k$. Because the
structure constants are rational numbers, Mal'cev Theorem
\cite{Malc} implies the existence of a discrete subgroup $\Gamma$
of $G$ such that the quotient space $M =\Gamma\backslash G$ is
compact. The cohomology of $M$ is given by
 \begin{eqnarray*}
 H^0(M) &=& \la 1\ra,\\
 H^1(M) &=& \la [\alpha_1], [\alpha_2],[\alpha_3]\ra,\\
 H^2(M) &=& \la [\alpha_{13}],
 [\alpha_{23}],[\alpha_{24}], [\alpha_{16} +\alpha_{25}- \alpha_{34}],
 [\alpha_{26}- \alpha_{45}] \ra,\\
  H^3(M) &=& \la [\alpha_{126}],
 [\alpha_{135}],[\alpha_{136}+\alpha_{146}], [\alpha_{136} +\alpha_{235}],
 [\alpha_{236}+ \alpha_{345}], [\alpha_{156}-\alpha_{246}+\alpha_{345}]
 \ra, \\
 H^4(M) &=& \la [\alpha_{2345}],[\alpha_{1236}],[\alpha_{2456}],
 [\alpha_{1456}+\alpha_{2346}],[\alpha_{1356}+\alpha_{1456}] \ra,\\
 H^5(M) &=& \la [\alpha_{23456}],[\alpha_{13456}],[\alpha_{12456}]\ra, \\
 H^6(M) &=& \la [\alpha_{123456}]\ra.
 \end{eqnarray*}

Therefore $M$ is a symplectic manifold with symplectic form
$\omega = \alpha_{16} + \alpha_{25} - \alpha_{34}$, and
$b_3(M)=6$. It is simple to check that $\Lo^2: H^1(M) \to H^5(M)$
is the zero map. On the other hand, $\Lo:H^1(M)\to H^3(M)$ has
kernel of dimension $1$ and generated by $[\alpha_1]$. This
follows from $\omega\wedge\alpha_1=d(\alpha_{45}+\alpha_{35})$, so
$[\alpha_1]$ is in the kernel, and
$[\omega\wedge\alpha_2\wedge\alpha_3]\neq 0$, so
$[\alpha_2],[\alpha_3]$ are not in the kernel. By Lemma
\ref{lem:hr3}, $b_3^{\rm hr}(M,\omega)= 6+1-3=4$.

But $M$ is of dimension $6$, and we need a manifold of dimension
$2n$, where $n-r>3$. We shall fix $2n=8+2r$ and define the
$2n$--dimensional manifold
 $$
 X=M\times\CP^{r+1}.
 $$
Let $\omega_0$ be the Fubini-Study symplectic form of $\CP^{r+1}$,
so $X$ is a sympletic manifold with symplectic form
$\Omega=\omega+\omega_0$. Now
  \begin{eqnarray*}
  H^1(X) &= & H^1(M) ,\\
  H^3(X) &= & H^3(M)\oplus (H^1(M)\otimes H^{2}(\CP^{r+1})),\\
  \vdots \\
  H^{2n-3}(X) &=&  (H^5(M)\otimes H^{2r}(\CP^{r+1})) \oplus  (H^3(M) \otimes
  H^{2r+2}(\CP^{r+1})) ,\\
  H^{2n-1}(X) &=&  H^5(M) \otimes H^{2r+2}(\CP^{r+1}).
 \end{eqnarray*}

First we will compute $b_3^{\rm hr}(X,\Omega)$ by using Lemma
\ref{lem:hr3}. Clearly $b_3(X)=6+3=9$. The map $\LO=\Lo +
L_{[\omega_0]}$, so $\LO^{n-1}: H^1(X) \to H^{2n-1}(X)=H^5(M)
\otimes H^{2r+2}(\CP^{r+1})$ equals
  $$
  \LO^{n-1} =\left(\Lo + L_{[\omega_0]}\right)^{n-1}=  \sum_{ j} {n-1 \choose j}
  \Lo^j L_{[\omega_0]}^{n-1-j}=0,
  $$
since $\Lo^j=0$ for $j>1$ and $L_{[\omega_0]}^{n-1-j}=0$ for
$n-1-j>r+1$, i.e.,\ for $j<2$. The map $\LO^{n-2}: H^1(X) \to
H^{2n-3}(X)= \left( H^5(M)\otimes H^{2r}(\CP^{r+1})\right) \oplus
\left( H^3(M) \otimes H^{2r+2}(\CP^{r+1})\right)$ equals
  $$
  \LO^{n-2} =  \sum_{ j} {n-2 \choose j}
  \Lo^j L_{[\omega_0]}^{n-2-j}=\Lo L_{[\omega_0]}^{r+1}.
  $$
So $\ker \left(\LO^{n-2}:H^1(X) \to H^{2n-3}(X)\right)= \ker
\left(\Lo:H^1(M) \to H^{3}(M)\right)=\la [\alpha_1]\ra$. Lemma
\ref{lem:hr3} yields
 $$
 b_3^{\rm hr}(X,\Omega)= 9+1-3 =7,
 $$
for any value of $r$. With these preliminaries at hand, we are
ready to start with our examples.

\begin{example}\label{example1}
The compact symplectic manifold $(X=M\times\CP^{2}, \Omega)$ has
an Auroux submanifold $Z_1\subset (X,\Omega)$ such that $b_3^{\rm
hr}(Z_1,\Omega_{Z_1}) > b_3^{\rm hr}(X,\Omega)$.
\end{example}

\begin{proof}
Let $A= [\alpha_{26}  - \alpha_{45}] \in H^2(M)$. To define an
Auroux submanifold $Z_1\subset (X,\Omega)$ in the conditions
required, we consider a rank $1$ bundle $E$ with first Chern class
$c_1(E) = A\in H^2(M)\subset H^2(X)$. Note that $n=5$, $r=1$ in
this case. Hence the Auroux submanifold $Z_1\subset X$ has
$PD[Z_1]= k[\Omega]+A$. To apply the part {\it (ii)} of Lemma
\ref{lem:hr3}, we need to compute the map $\LO^{3}(k \LO + L_A):
H^1(X) \to H^{9}(X)=H^5(M) \otimes H^{4}(\CP^{2})$, where $L_A$ is
the map in cohomology given by cup product with the class $A$.
This is
  $$
  \LO^{3}(k \LO + L_A) =\LO^{3}L_A =  \Lo L_A  L_{[\omega_0]}^{2},
  $$
since $\LO^{4}=0$, by the above calculation. This map has kernel
of dimension $1$, generated by $[\alpha_1]$, since $\Lo
([\alpha_1])=0$, but
 \begin{equation}\label{eqn:A}
 [\alpha_2]\cup [\alpha_3]\cup [\omega]\cup A \neq 0.
 \end{equation}
The map $\LO^{2}(k\LO+A): H^1(X) \to H^{7}(X)=(H^5(M)\otimes
H^{2}(\CP^{2}) )\oplus (H^3(M) \otimes H^{4}(\CP^{2}))$ equals
  $$
  \LO^{2}(k\LO+L_A) =  k \Lo L_{[\omega_0]}^{2}  + \Lo L_A L_{[\omega_0]}.
  $$
The first component has kernel generated by $[\alpha_1]$, by what
we have seen above. The second component has the same kernel
again, so $\dim \ker \left(\LO^{2}(k\LO+L_A): H^1(X) \to
H^{7}(X)\right) =1$. Now Lemma \ref{lem:hr3} gives
 $$
 b_3^{\rm hr}(Z_1,\Omega_{Z_1})= 9+1-1 =9.
 $$
Therefore,  $b_3^{\rm hr}(Z_1,\Omega_{Z_1})> b_3^{\rm
hr}(X,\Omega)$.
\end{proof}

\bigskip
Notice that in the example above, all the calculation hinges in
(\ref{eqn:A}). In fact, we have

\begin{example}\label{example2}
The compact symplectic manifold $(X=M\times\CP^{2}, \Omega)$ has
an Auroux submanifold $Z_1'\subset (X,\Omega)$ such that $b_3^{\rm
hr}(Z_1',\Omega_{Z_1'})=b_3^{\rm hr}(X,\Omega)$.
\end{example}

\begin{proof}
We take a class $A\in H^2(M)$ such that $[\alpha_2]\cup
[\alpha_3]\cup [\omega]\cup A = 0$; for instance, use
$A=[\alpha_{13}]$. Then we obtain an Auroux submanifold
$(Z_1',\Omega_{Z_1'})$ of $(X,\Omega)$ with $b_3^{\rm
hr}(Z_1',\Omega_{Z_1'})= b_3^{\rm hr}(X,\Omega)=7$.
\end{proof}

\bigskip
Finally we give an example where the Auroux submanifold has less
harmonic cohomology than the ambient submanifold.

\begin{example}\label{example3}
There are Auroux submanifolds $(Z_3,\Omega_{Z_3})\subset
(Z_1,\Omega_{Z_1})\subset (X=M\times\CP^{4}, \Omega)$ such that
$b_3^{\rm hr}(Z_3,\Omega_{Z_3})< b_3^{\rm hr}(Z_1,\Omega_{Z_1})$.
\end{example}

\begin{proof}
Consider the manifold $(X=M\times\CP^{4}, \Omega)$ of dimension
$14$ (now $n=7$ and $r=3$). Take again $A= [\alpha_{26}  -
\alpha_{45}] \in H^2(M)$ and let $E$ be a rank $1$ bundle with
$c_1(E)=A$. There is another bundle $F$ such that $E\oplus F$ is a
trivial bundle. Actually, one may take $F$ to have rank $2$ and
Chern classes $c_1(F)=-A$ and $c_2(F)=A^2$. Let $Z_3\subset
(X=M\times\CP^{4}, \Omega)$ be the Auroux submanifold (of
codimension $6$) associated to the (trivial) bundle $E\oplus F$.
Since the Chern classes of $E\oplus F$ are all zero, we have that
$$
 b_3^{\rm hr}(Z_3,\Omega_{Z_3})= b_3^{\rm
 hr}(X,\Omega)=7,
$$
as in Example~\ref{example2}.

Now let $Z_1\subset (X=M\times\CP^{4}, \Omega)$ be the Auroux
submanifold associated to the bundle $E$. By
Example~\ref{example1}, we have that $b_3^{\rm
hr}(Z_1,\Omega_{Z_1})=9$. But, the construction in \cite{Au97} is
carried out in such a way that $Z_3$ is also an Auroux submanifold
of $Z_1$ (of codimension $4$), and
$$
 b_3^{\rm hr}(Z_3,\Omega_{Z_3})< b_3^{\rm
 hr}(Z_1,\Omega_{Z_1}).
$$
\end{proof}

\section{Symplectic blow-ups} \label{symplecticblowups}

This section is devoted to the study of the $s$--Lefschetz
property for the symplectic blow-up ${\wcpn}$ of the complex
projective space ${CP}^m$ along a symplectic submanifold $M \inc
{CP}^m$.

Let $(M,\omega)$ be a {\em compact\/} symplectic manifold of
dimension $2n$. Without loss of generality we can assume that the
symplectic form $\omega$ is integral (by perturbing it to make it
rational and then rescaling), i.e.,  $[\omega]\in H^2(M;\ZZ)$.
 A theorem of Gromov and Tischler
\cite{Gr1,Ti} states that there is a symplectic embedding $i\colon
(M,\omega)\longrightarrow (CP^m,\omega_0)$, with $m\geq {2n+1}$,
where $\omega_0$ is the standard K\"ahler form on ${CP}^m$ defined
by its natural complex structure and the Fubini--Study metric. We
take the symplectic blow-up $\wcpn$ of ${CP}^m$ along the
embedding $i$ (see \cite{McD1}). Then $\wcpn$ is a simply
connected compact symplectic manifold.

Recall that $i^*\omega_0=\omega$. We will denote also by
$\omega_0$ the pull back of $\omega_0$ to $\wcpn$ under the
natural projection $\wcpn\to CP^m$. Let $\widetilde{M}$ be the
projectivization of the normal bundle of the embedding $M\inc
CP^m$. Then $\pi\colon \widetilde{M}\longrightarrow M$ is a
locally trivial bundle with fiber $CP^{m-n-1}$. We will denote by
$\nu$ the Thom form of the submanifold $\widetilde{M}\subset
\wcpn$. The class $[\nu]$ is called the Thom class of the blow-up.
Then $\wcpn$ has a symplectic form $\Omega$ whose cohomology class
is $[\Omega]=[\omega_0]+\epsilon\, [\nu]$ for $\epsilon>0$ small
enough.

Let us consider a closed tubular neighborhood $\widetilde{W}$ of
$\widetilde{M}$ in $\wcpn$. By the tubular neighborhood theorem we
know that the normal bundle of $\widetilde{M}\inc \wcpn$ contains
a disk subbundle which is diffeomorphic to $\widetilde{W}$. Denote
by $p\colon {\widetilde{W}}\longrightarrow \widetilde{M}$ the
natural map. There is a map $q:\Omega^*(M)\to \Omega^{*+2}(\wcpn)$
given by pull-back by $\pi: \widetilde{M} \to M$, followed by
extending to a neighborhood of $\widetilde{M}$ using $p:
\widetilde{W}\to \widetilde{M}$ and then wedging by $\nu$, i.e.,\
$q(\alpha)=p^*\pi^*(\alpha)\wedge \nu$. We shall denote
$q(\alpha)=\alpha\wedge\nu$ for short. Note that
 $$
  (\alpha\wedge\nu) \wedge( \beta\wedge\nu) =
  (\alpha\wedge\beta\wedge \nu)\wedge \nu\, ,
 $$
for $\alpha,\beta\in\Omega^*(M)$. This makes notations of the type
$\alpha\wedge\beta\wedge\nu^2$ unambiguous. Also remark that
$[\omega_0\wedge\nu]=[\omega\wedge \nu]$ although
$\omega_0\wedge\nu\neq \omega\wedge\nu$ as forms.

The cohomology of $\wcpn$ was studied by McDuff~\cite{McD1}. There
she proved that there is a short exact sequence
\begin{equation} \label{eqn:v7}
  0\longrightarrow H^*(CP^m)\longrightarrow H^*(\wcpn)
   \longrightarrow A^*\longrightarrow 0,
\end{equation}
where $A^*$ is a free module over $H^*(M)$ generated by $\{[\nu],
[\nu^2],\cdots,[\nu^{m-n-1}]\}$.

Before going on to the study of the $s$--Lefschetz property for
$\wcpn$, we need to recall  the splitting of the cohomology groups
in terms of the primitive classes proved by Yan \cite{Ya} for hard
Lefschetz symplectic manifolds. His proof also works for
$s$--Lefschetz symplectic manifolds.

\begin{lemma} \label{primitiveclasses}
Let $(M,\omega)$ be a compact symplectic manifold of dimension
$2n$ satisfying the $s$--Lefschetz property for $s\leq n-1$. Then,
there is a splitting
 $$
 H^k(M)=P_k(M)\oplus L(H^{k-2}(M)),
 $$
where $P_k(M)$ is given by
 $$
 P_k(M) =\{v\in H^k(M) \ \mid \  L^{n-k+1}(v)=0\},
 $$
for $k\leq s$. The elements in $P_k(M)$ are called primitive
cohomology classes of degree $k$.
\end{lemma}

\begin{proof}
First, let us see that $P_k(M)\cap \Im L=0$. Take $x\in P_k(M)$
with $x=L(y)$, $y\in H^{k-2}(M)$. Then
$L^{n-k+2}(y)=L^{n-k+1}(x)=0$. By the $(k-2)$--Lefschetz property,
$y=0$ and hence $x=0$.

Now let us consider $a\in H^k(M)$ with $k\leq s$, and take the
element $L^{n-k+1}(a)\in H^{2n-k+2}(M)$. If $L^{n-k+1}(a)$ is the
zero class, then $a\in P_k(M)$ and the lemma is proved. If
$L^{n-k+1}(a)$ is non-zero, then there exists $b\in H^{k-2}(M)$
such that $L^{n-k+1}(a)=L^{n-k+2}(b)$ since $(M,\omega)$ is
$s$--Lefschetz and so the map $L^{n-k+2}: H^{k-2}(M)
\longrightarrow H^{2n-k+2}(M)$ is an isomorphism. Hence
${a-L(b)}\in P_k(M)$. But $a=(a-L(b))+L(b)$ which lies in
$P_k(M)\oplus \Im L$.
\end{proof}

According Lemma \ref{primitiveclasses} we can write
 \begin{equation}\label{eqn:cc2}
 H^k(M)=P_k(M)\oplus ( P_{k-2}(M)\cup [\omega]) \oplus
 \cdots \oplus (P_{k-2\lambda}(M)\cup [\omega^{\lambda}]),
 \end{equation}
with $\lambda = [\frac{k}{2}]$.

\begin{theorem} \label{lefschetz-blowup}
For any $s\leq n-1$, if $(M,\omega)$ is $s$--Lefschetz then there
exists $\epsilon_0>0$ such that $(\wcpn,\Omega=\omega_0+\epsilon
\nu)$ is $(s+2)$--Lefschetz, for any $\epsilon \in
(0,\epsilon_0]$. In particular, for $\epsilon\in \QQ\cap
(0,\epsilon_0]$, we have that $[\Omega_\epsilon]$ a rational class
(and hence a multiple of it is integral).
\end{theorem}

\begin{proof}
Following the notation stated at the beginning of this Section, we
must prove that the map $[\omega_0+\epsilon \nu]^{m-k} \colon
H^k(\wcpn) \longrightarrow H^{2m-k}(\wcpn)$ is an isomorphism  for
any $k\leq s+2 \leq n+1$. First, using (\ref{eqn:v7}) and
(\ref{eqn:cc2}), we notice that for $k\leq s+2$ the cohomology
group $H^k(\wcpn)$ is generated by the classes:
 $$
 \left\{\begin{array}{ll} [\omega_0]^{\frac {k}{2}}, & \mbox{if } k \mbox{ is
 even},\\[12pt]
 [p_{k-2i-2t} \wedge \omega_0^t \wedge \nu^i], & \mbox{where }
 [p_{k-2i-2t}]\in P_ {k-2i-2t}(M),
   i>0, t\geq 0  \mbox{ and }  i+t \leq [\frac{k}{2}].
 \end{array}\right.
 $$
Suppose that $k$ is even (the proof is similar when $k$ is odd).
We prove that the map $[\omega_0+\epsilon \nu]^{m-k}$ is injective
by computing each one of the following cohomology classes in
$H^{2m}(\wcpn)$:
 $[\omega_0+\epsilon \nu]^{m-k} \cup
 [\omega_0]^{\frac {k}{2}} \cup [\omega_0]^{\frac {k}{2}} $,
 $[\omega_0+\epsilon \nu]^{m-k} \cup  [p_{k-2i-2t}\wedge \omega_0^t \wedge \nu^i]
  \cup [\omega_0]^{\frac {k}{2}}$ for $i+t \leq {\frac{k}{2}}$, and
 $[\omega_0+\epsilon \nu]^{m-k} \cup  [p_{k-2i-2t}\wedge \omega_0^t \wedge \nu^i]
 \cup  [q_{k-2j-2s}\wedge \omega_0^s \wedge \nu^j]$ if  $i+t, j+s \leq {\frac{k}{2}}$,
 where $[q_{k-2j-2s}]\in P_ {k-2j-2s}(M)$.

We begin by showing that the class $[\omega_0+\epsilon \nu]^{m-k}
\cup [\omega_0]^{\frac {k}{2}} \cup [\omega_0]^{\frac {k}{2}} $ is
non-trivial. We have
 $$
 [\omega_0+\epsilon \nu]^{m-k} \cup [\omega_0]^{\frac {k}{2}} \cup [\omega_0]^{\frac {k}{2}} =
 \sum_{r= 0}^{m-k} {m-k \choose r} \epsilon^r [\omega_0^{m-r} \wedge \nu^r ]=
 [\omega_0]^{m} + \sum_{r= 1}^{m-k} {m-k \choose r} \epsilon^r [\omega_0^{m-r} \wedge \nu^r].
 $$
In this sum, the terms $[\omega_0^{m-r} \wedge \nu^r]$ are zero
for $1\leq r \leq m-n-1$ since $M$ has  dimension $2n$ and so
$[\omega_0^{n+1}\wedge \nu]= [\omega^{n+1}\wedge \nu]=0$. Then,
 \begin{eqnarray} \label{eqn:aa}
 [\omega_0+\epsilon \nu]^{m-k} \cup [\omega_0]^{\frac {k}{2}} \cup [\omega_0]^{\frac {k}{2}}& =&
 [\omega_0]^{m} + {m-k \choose {m-n}} \epsilon^{m-n} [\omega_0^{n} \wedge \nu^{m-n}] \nonumber\\ &  &
 +\sum_{r= m-n+1}^{m-k} {m-k \choose r} \epsilon^r [\omega_0^{m-r} \wedge \nu^r]
 \\
 & =&
  [\omega_0]^{m} + {m-k \choose {m-n}} \epsilon^{m-n} [\omega^n\wedge \nu^{m-n}]
  + O(\epsilon^{m-n+1}), \nonumber
 \end{eqnarray}
which is a non-zero class (for $\epsilon$ small enough).

Proceeding in a similar way, let $i+t \leq \frac{k}{2}$, $i>0$,
$t\geq 0$, and $[p_{k-2i-2t}]\in P_{k-2i-2t}(M)$. Then
   \begin{eqnarray}\label{eqn:bb}
 [\omega_0+\epsilon \nu]^{m-k} \cup [p_{k-2i-2t}\wedge \omega_0^t \wedge \nu^i]
  \cup [\omega_0]^{\frac {k}{2}}
 & =&
   \sum_{r= 0}^{m-k} {m-k \choose r} \epsilon^r  [p_{k-2i-2t} \wedge
   \omega_0^{t+m-\frac{k}{2}-r}\wedge \nu^{r+i}] \nonumber\\
 &=&  {m-k \choose {m-n-i}} \epsilon^{m-n-i}
 [p_{k-2i-2t}\wedge \omega^{n+ i+ t -{\frac {k}{2}}}  \wedge
 \nu^{m-n}] \\[12pt] & &
 + \, O(\epsilon^{m-n-i+1}), \nonumber
 \end{eqnarray}
using that for $i < m-n-r$, we have that $[p_{k-2i-2t} \wedge
\omega^{t+m-\frac{k}{2}-r}\wedge\nu^{r+i}]=0$, since $\deg
(p_{k-2i-2t} \wedge \omega^{t+m-\frac{k}{2}-r})>2n$. Suppose that
 \begin{equation} \label{eqn:bb2}
 x=a[\omega_0]^{\frac{k}{2}} + \sum_{i+t\leq \frac{k}{2}, i>0}
 [p_{k-2i-2t}\wedge \omega_0^t \wedge \nu^i]
 \in H^{k}(\wcpn)
 \end{equation}
is an element such that $ [\omega_0+\epsilon \nu]^{m-k} \cup x
=0$. Then multiplying by $[\omega_0]^{\frac{k}{2}}$ and using
(\ref{eqn:aa}) and (\ref{eqn:bb}), we get that $a=0$. So
 \smallskip
 \begin{equation} \label{eqn:cc}
  x= \sum_{i+t\leq \frac{k}{2}, i>0}
  [p_{k-2i-2t}\wedge \omega_0^t \wedge \nu^i]\, .
 \end{equation}

Now we compute for $i+t\leq \frac{k}{2}$ and $j+s\leq \frac{k}{2}$
the following product
 \begin{eqnarray}\label{eqn:dd}
  &[\omega_0+\epsilon \nu]^{m-k} \cup  [p_{k-2i-2t}\wedge \omega_0^t \wedge \nu^i]
 \cup  [q_{k-2j-2s}\wedge \omega_0^s \wedge \nu^j]
  =  \\[12pt]
  &=\displaystyle {m-k \choose {m-n-i-j}} \epsilon^{m-n-i-j}  [p_{k-2i-2t} \wedge q_{k-2j-2s}
   \wedge \omega^{n-k+i+t+j+s} \wedge\nu^{m-n} ]
  +O(\epsilon^{m-n-i-j+1}). \nonumber
 \end{eqnarray}
Let us concentrate on the leading term. The duality on $H^r(M)$
defines a duality on the space $P_r(M)$ of the primitive
cohomology classes:
 $$
 p^{\sharp} \colon P_r(M)\otimes P_r(M) \longrightarrow \ {\RR}
 $$
given by
 $$
 p^{\sharp} ([\alpha],[\beta])=\int_M \alpha\wedge\beta\wedge\omega^{n-r},
 $$
which is nondegenerate, but
 $$
 p^{\sharp} \colon P_r(M)\otimes P_{r+2s}(M) \longrightarrow \ {\RR}
 $$
given by
 $$
 p^{\sharp} ([\alpha],[\beta])=\int_M \alpha\wedge\beta\wedge\omega^{n-r-s},
 $$
is zero if $s\neq 0$, since $[\omega]^{n-r-s}$ maps $P_{r+2s}(M)$
to zero. Thus the matrix $A_{i+t,j+s}$ associated to
$p^{\sharp}:P_{k-2i-2t}(M)\otimes P_{k-2j-2s}(M)\to \RR$ is
non-singular if $i+t=j+s$ and zero if $i+t\neq j+s$.

Consider the spaces
 $$
  P_\mu := \bigoplus_{i+t=\mu, i>0} P_{k-2i-2t}(M)\,
  [\omega^t]\,[\nu^i]
 $$
and
 $$
  W=\bigoplus_{1\leq \mu\leq \frac{k}{2}} P_\mu\, ,
  $$
so that $H^k(\wcpn)= [\omega_0^{\frac{k}{2}}] \oplus W$. There is
a bilinear map
 $$
 p^{\sharp}_1 \colon W\otimes W \longrightarrow \ {\RR}
 $$
given by
 $$
 p^{\sharp}_1  ([p_{k-2i-2t}\wedge \omega_0^t \wedge \nu^i],
 [q_{k-2j-2s}\wedge \omega_0^s \wedge \nu^j]) = \int_{\wcpn}
  p_{k-2i-2t}\wedge  q_{k-2j-2s}\wedge \omega^{n-k+i+t+j+s} \wedge
  \nu^{m-n}\, .
 $$
The matrix $B_\mu$ of $p^{\sharp}_1|_{P_\mu\otimes P_\mu}$ is the
block matrix whose block in the place $(i,j)$ with $1\leq i,j\leq
\mu$ is the matrix
 $$
 {m-k \choose {m-n-i-j}}\cdot \epsilon^{m-n-i-j} \cdot A_\mu \, .
 $$
Let $d=\dim P_{k-2i-2t}(M)$. The determinant  of $B_\mu$  is
 \begin{eqnarray} \label{eqn:formul}
 \det(A_\mu)^{\mu} \cdot \left[ \det \left( \epsilon^{m-n-i-j} {m-k \choose
 {m-n-i-j}}\right)_{1\leq i,j\leq \mu} \right]^d = \qquad
 \\ \qquad =
 \det(A_\mu)^{\mu} \cdot
 \left[\epsilon^{(m-n)\mu-\mu(\mu+1)}
 \frac{\displaystyle {{m-k+\mu-1}\choose{m-n-\mu-1}} \cdots
 {{m-k}\choose{m-n-\mu-1}}}
   {\displaystyle {{m-n-2}\choose{m-n-\mu-1}}\cdots{{m-n-\mu-1}\choose{m-n-\mu-1}}
   }\right]^d, \nonumber
 \end{eqnarray}
which is of the form $\lambda_\mu\cdot\epsilon^{a_\mu}$ where
$\lambda_\mu\neq 0$. Here we use that $k\leq s+2\leq n+1
\Rightarrow m-k > m-n-\mu-1$ and $\mu \leq \frac{k}{2} < m-n
\Rightarrow m- n-\mu -1\geq 0$.

The determinant of the matrix of $p^{\sharp}_1$ is the product of
$\det B_\mu$ for $1\leq \mu\leq \frac{k}{2}$, hence of the form
$\lambda \cdot\epsilon^{a}$ where $\lambda\neq 0$. The matrix
associated to the bilinear map $p^{\sharp}_2 \colon W\otimes W
\longrightarrow \ {\RR}$ given by
 $$
 p^{\sharp}_2  ([p_{k-2i-2t}\wedge \omega_0^t \wedge \nu^i],
 [q_{k-2j-2s}\wedge \omega_0^s \wedge \nu^j]) =
 [\omega_0+\epsilon \nu]^{m-k} \cup [p_{k-2i-2t}\wedge \omega_0^t
 \wedge \nu^i] \cup [q_{k-2j-2s}\wedge \omega_0^s \wedge \nu^j]
 $$
has at each entry an $\epsilon$--perturbation of the corresponding
entry of $B_\mu$, by (\ref{eqn:dd}). Hence its determinant is
$\lambda\cdot\epsilon^a + O(\epsilon^{a+1})$ and it is nonzero for
small $\epsilon>0$. Therefore $p^{\sharp}_2$ is a pairing and
hence (\ref{eqn:cc}) is zero. So $\wcpn$ is $(s+2)$--Lefschetz.

To complete the proof, we must notice that in the conditions of
Theorem \ref{lefschetz-blowup}, there exists $\epsilon_0>0$ such
that for any $\epsilon\in (0,\epsilon_0]$ the manifold
$(\widetilde{CP}^m,\Omega_\epsilon=\omega_0 + \epsilon \nu)$ is
$(s+2)$-Lefschetz. In particular, if $[\omega_0]$ is an integral
$2$-cohomology class, then for rational $\epsilon>0$, we have that
$[\Omega_\epsilon]$ is a rational class, hence a multiple of it is
an integral class.
\end{proof}

\begin{remark}
Cavalcanti \cite[Theorem 4.2]{Cav} has proved that if $M$ is hard
Lefschetz then $\wcpn$ is also hard Lefschetz. This also can be
proved with the arguments of Theorem~\ref{lefschetz-blowup} with
few modifications:

We suppose $M$ is hard-Lefschetz and must prove that $\wcpn$ is
$k$--Lefschetz for any $n+2\leq k\leq m-1$. In this case, the
group $H^k(\wcpn)$ is generated by $[\omega_0]^{\frac {k}{2}}$ (if
$k$ is even) and $[p_{k-2i-2t} \wedge \omega_0^t \wedge \nu^i]$,
$[p_{k-2i-2t}]\in P_{k-2i-2t}(M)$, $0<i<m-n$, $k-n \leq t+2i$,
$t+i \leq [\frac{k}{2}]$. The rest of the argument is unchanged
except at two points: use that $i<m-n$ in (\ref{eqn:bb}) to get
that $a=0$ in (\ref{eqn:bb2}); and use that $2\mu \geq k-n
\Rightarrow m-k \geq m-n-\mu-1$ to get that $\lambda_\mu \neq 0$
in (\ref{eqn:formul}).
\end{remark}

The following result shows that the converse of the previous
theorem is also true if $M$ is parallelizable.

\begin{proposition} \label{nolefschetz-blowup}
Let $(M,\omega)$ be a compact symplectic manifold of dimension
$2n$, such that $M$ is parallelizable and $(M,\omega)$ is not
$s$--Lefschetz for some $s\geq 1$. Then $\wcpn$ is not
$(s+2)$--Lefschetz.
\end{proposition}

\begin{proof}
Since $M$ is parallelizable, its tangent bundle $TM$  is trivial.
Denote by $N$ the normal bunble of $M\inc CP^m$. Then the
restriction to $M$ of the tangent bundle of $CP^m$  is
$TCP^m\vert_M=TM \oplus N$. The total Chern class of $N$ is given
by $c(N)=c(TCP^m\vert_M)= (1+[\omega])^{m+1}$, so $c_i(N)$ is a
multiple of $[\omega]^i$.

Taking into account  that $(M,\omega)$ is not $s$--Lefschetz, we
know that there is a non-trivial class $[p_s]\in H^s(M)$ such that
$[p_s]\in \ker (H^s(M) \times H^s(M) \longrightarrow \ {\RR})$.
This means that for any other element $[q_s] \in H^s(M)$ we have
that $[p_s\wedge q_s\wedge \omega^{n-s}]=0$ in $H^*(M)$. In the
cohomology ring $H^*(\wcpn)$ we have the following equality
 $$
 [ p_s\wedge \nu \wedge q_s \wedge \omega^l_0\wedge \nu^{m-s-l-1}]=
  \left\{ \begin{array}{ll}  0,  & \hbox{if } m-s-l < m-n,  \\
   {} [p_s\wedge q_s\wedge\omega^{n-s} \wedge \nu^{m-n}]=0, \qquad &
   \hbox{if } m-s-l=m-n,
  \\
  {} [p_s\wedge q_s\wedge\omega^l\wedge P(c(N)) \wedge
  \nu^{m-n}]=0, & \hbox{if } m-s-l >m-n, \end{array} \right.
  $$
since $P(c(N))$ is a polynomial in the Chern classes of $N$, and
hence a multiple of $[\omega]^{n-s-l}$, because the Chern classes
of $N$ are multiples of powers of $[\omega]$.

Therefore for any $j+l\leq \frac{s+2}{2}$, $j>0$, and
$[q_{s+2-2j-2l}\wedge \omega_0^l \wedge \nu^j] \in
H^{s+2}(\wcpn)$, we have
 $$
 [\omega_0+\epsilon \nu]^{m-s-2} \cup [p_s \wedge \nu]\cup
 [q_{s+2-2j-2l} \wedge \omega_0^l \wedge \nu^j]=0.
 $$
Also, in the case where $s+2$ is even, we have
 $$
 [\omega_0+\epsilon \nu]^{m-s-2} \cup [p_s \wedge \nu] \cup
 [\omega_0]^{\frac{s+2}{2}}=0.
 $$

Thus $[p_s\wedge \nu]\in \ker (H^{s+2}(\wcpn)\times H^{s+2}(\wcpn)
\longrightarrow \ {\RR})$, which  proves that $\wcpn$ is not
$(s+2)$--Lefschetz.
\end{proof}

\section{Examples of $s$-Lefschetz symplectic manifolds} \label{Donaldsonhr}

In this section, examples of compact symplectic manifolds which
are $s$--Lefschetz but not $(s+1)$--Lefschetz are constructed for
$s=3$ and for any {\em even} integer $s\geq 2$.

First we show the existence of a simply connected compact
symplectic manifold $M_s$, of high dimension,  which is
$s$--Lefschetz but not $(s+1)$--Lefschetz, for each even integer
value of $s\geq 2$. The idea for the construction of $M_s$ is to
follow an iterative procedure starting from an appropriate low
dimensional compact symplectic manifold, take a symplectic
embedding of it in a complex projective space $CP^{m}$ and then
consider the symplectic blow-up of $CP^{m}$ along the embedded
submanifold in order to get a simply connected compact symplectic
manifold which, according to Theorem~\ref{lefschetz-blowup}, will
be Lefschetz up to a strictly higher level.

The starting point to construct $M_s$ will be the
Kodaira--Thurston manifold $KT$~\cite{Kod, Tu}. We begin reviewing
it. Consider the Heisenberg group $H$, that is, the connected
nilpotent Lie group of dimension $3$ consisting of matrices of the
form
 $$
 a=\left( \begin{array}{ccc} 1&x&z\\ 0&1&y\\ 0&0&1 \end{array} \right) ,
 $$
where $x,y,z \in {\RR}$. A global system of coordinates $(x,y,z)$
for $H$ is given by $x(a)=x$, $y(a)=y$, $z(a)=z$, and a standard
calculation shows that $\{dx,\ dy,\ dz-xdy\}$ is a basis for the
left invariant $1$--forms on $H$. Let $\Gamma$ be the discrete
subgroup of $H$ consisting of matrices whose entries $x$, $y$ and
$z$ are integer numbers. So the quotient space $\Gamma{\backslash}
H$ is compact, and the forms $dx$, $dy$, $dz-xdy$ descend to
$1$--forms $\alpha$, $\beta$, $\gamma$ on $\Gamma{\backslash} H$
such that $\alpha$ and $\beta$ are closed, and $d\gamma=-\alpha
\wedge \beta$.

The {\it Kodaira--Thurston manifold} $KT$ is the product $KT
=\Gamma{\backslash} H \times S^1$ (see~\cite{Kod, Tu}). Now, if
$\eta$ is the standard invariant $1$--form on $S^1$, then
$\{\alpha, \beta, \gamma, \eta\}$ constitutes a (global) basis for
the $1$--forms on $KT$. Since
 $$
 d\alpha=d\beta=d\eta=0, \quad  d\gamma=-\alpha \wedge \beta ,
 $$
using Nomizu's theorem~\cite{No} we compute the real cohomology of
$KT$:
\begin{eqnarray*}
 H^0(KT) &=& \la 1\ra, \\
 H^1(KT) &=& \la [\alpha], [\beta], [\eta]\ra,\\
 H^2(KT) &=& \la [\alpha \wedge \gamma], [\beta\wedge \gamma],
   [\alpha \wedge \eta], [\beta \wedge \eta]\ra,\\
 H^3(KT) &=& \la [\alpha \wedge \gamma \wedge \eta],
   [\beta\wedge \gamma \wedge \eta], [\alpha \wedge \beta \wedge \gamma]\ra,\\
 H^4(KT) &=& \la [\alpha \wedge \beta \wedge \gamma \wedge \eta]\ra.
\end{eqnarray*}

Therefore, $KT$ is a symplectic manifold with the symplectic form
$\omega= \alpha \wedge \gamma + \beta \wedge \eta$. It is clear
that $(KT,\omega)$ is not $1$--Lefschetz, which follows directly
from its cohomology or from the general result of Benson and
Gordon~\cite{BG}. Moreover, $H^k_{\rm hr}(KT,\omega)=H^k(KT)$ for
any $k\not= 3$, but $b_3^{\rm hr}(KT,\omega)=2<3=b_3(KT)$. It is
easy to see that the same holds for any other symplectic form on
$KT$.

Denote $M_0=KT$. By Gromov--Tischler theorem~\cite{Gr1, Ti} there
exists a symplectic embedding of $(KT,\omega)$ in the complex
projective space $CP^{m_0}$, with $m_0=5$, endowed with its
standard K\"ahler form. Let us denote by
$(M_2=\widetilde{CP}{}^{m_0},\Omega_2)$ the blow-up of $CP^{m_0}$
along $M_0$. By Theorem~\ref{lefschetz-blowup} we can consider
$\Omega_2$ an integral form. We may again embed symplectically
$(M_2,\Omega_2)$ into $CP^{m_2}$ with $m_2=11$ and blow-up
$CP^{m_2}$ along $M_2$ to obtain
$(M_4=\widetilde{CP}{}^{m_2},\Omega_4)$. So in this fashion we get
a simply connected compact symplectic manifold $(M_s,\Omega_s)$
for any even integer $s\geq 2$ obtained as the symplectic blow-up
$\widetilde{CP}{}^{m_{s-2}}$ of $CP^{m_{s-2}}$ along
$(M_{s-2},\Omega_{s-2})$ symplectically embedded into
$CP^{m_{s-2}}$, where $m_{s-2}=2m_{s-4}+1$. Notice that the
dimension of the manifold $M_{s+2}$ is equal to $2m_s$, where
 $$
 m_s= 6\cdot 2^r - 1,
 $$
for $s=2r\geq 0$.

\begin{proposition} \label{2slefschetz-Ms}
For any even integer $s\geq 2$, the simply connected compact
symplectic manifold $M_s=\widetilde{CP}{}^{m_{s-2}}$ is
$s$--Lefschetz but not $(s+1)$--Lefschetz.
\end{proposition}

\begin{proof}
Since $M_0=KT$ is $0$--Lefschetz (any symplectic manifold is), we
can apply Theorem~\ref{lefschetz-blowup}  $r$ times, with $2r=s$,
to conclude that the manifold $M_{s}$ is $s$--Lefschetz. To show
that $M_{s}$ is not $(s+1)$--Lefschetz we note the following fact.
Consider $(M,\omega)$ a compact symplectic manifold and embed
symplectically $M \inc CP^{m}$ with $m \geq 2n+1$, where $2n$ is
the dimension of $M$. As  usual we write $\wcpn$ for the
symplectic blow-up of $CP^{m}$ along $M$. By (\ref{eqn:v7}), the
Betti number $b_i(\wcpn)$ is given by
 $$
 b_i(\wcpn)=b_{i-2}(M)+b_{i-4}(M)+\cdots+b_1(M)
 $$
if $i>1$ is odd. Therefore, $b_3(M_2)=b_1(KT)=3$. For $M_4$, we
have $b_1(M_4)=b_3(M_4)=0$ and $b_5(M_4)=3$. In general, for any
manifold $M_s$ the odd Betti numbers $b_{2j-1}(M_s)$ vanish for
$j\leq r$, and $b_{s+1}(M_s)=b_1(KT)=3$. This proves that $M_s$ is
not $(s+1)$--Lefschetz using Proposition~\ref{Bettinumbers}.
\end{proof}

In the following result we decrease as much as possible the
dimension of the examples constructed in
Proposition~\ref{2slefschetz-Ms} by using iterated Donaldson
symplectic submanifolds.

\begin{proposition}\label{2slefschetz}
Let $s\geq 2$ be an even integer, and let $M_s$ be the simply
connected compact symplectic manifold constructed in
Proposition~$\ref{2slefschetz-Ms}$. Then, there is a symplectic
submanifold $W_s \inc M_s$ of dimension $2(s+2)$ which is
$s$--Lefschetz but not $(s+1)$--Lefschetz, and every de Rham
cohomology class in $H^i(W_s)$ admits a symplectically harmonic
representative for any~$i\not= s+3$.
\end{proposition}

\begin{proof}
According to Theorem \ref{lefschetz-blowup}, we can assume that
the symplectic form $\Omega_s$ of $M_s$ is an integral form and
$(M_s,\Omega_s)$ is $s$--Lefschetz. Therefore, we can consider an
iterated Donaldson symplectic submanifold $Z_l\inc M_s$ of
codimension $2l$, i.e. $\dim Z_l=2(m_{s-2}-l)$. In particular, if
$s=2r$ then we take $l_s=m_{s-2}-s-2= 6\cdot 2^{r-1} - 2r -3$, and
denote by $W_s$ the corresponding simply connected compact
symplectic manifold $Z_{l_s}$ of dimension $2(s+2)$.

Since $6\cdot 2^r -2r -3=2m_{s-2}-s-1$, Poincar\'e duality implies
that $b_{6\cdot 2^r -2r -3}(M_s)=b_{s+1}(M_s)$, which equals
$b_1(KT)=3$ as shown in the proof of
Proposition~\ref{2slefschetz-Ms}.

Notice that $6\cdot 2^r -2r -3 = s+3+2l_s$. Therefore,
$b_{s+3}(W_s)=b_{s+3+2l_s}(M_s)=3$. Moreover,
Corollary~\ref{q-submanifold-2} implies that $b_i(W_s) - b_i^{\rm
hr}(W_s)=0$ for $i> (s+3)$, and $b_{s+3}(W_s) - b_{s+3}^{\rm
hr}(W_s) = b_{s+3+2l_s}(M_s) - b_{s+3+2l_s}^{\rm hr}(M_s) \equiv 1
\pmod 2$, by Proposition \ref{Bettinumbers}. {}From
Proposition~\ref{grados+2} we conclude that $W_s$ is
$s$--Lefschetz but not $(s+1)$--Lefschetz.
\end{proof}

\begin{remark}
If we begin with {\em any} symplectic $4$--manifold $N$ whose
first Betti number is $b_1(N)=1$ (see~{\rm \cite{Go}}), then we
obtain a symplectic manifold $W'_s$ satisfying the conditions of
Proposition~\ref{2slefschetz}, but with $b_{s+3}^{\rm
hr}(W'_s)=0$.
\end{remark}

\begin{corollary}\label{high-dim} Let $n$ and $s$ be
integer numbers such that $s\geq 2$ is even, and $n\geq s+2$. Then
there exists a simply connected compact symplectic manifold of
dimension~$2n$ which is $s$--Lefschetz but not $(s+1)$--Lefschetz.
\end{corollary}

It is worthy to remark that Proposition~\ref{2slefschetz} and
Corollary~\ref{high-dim} also hold in the {\em non-simply
connected\/} setting. For any even integer $s\geq 2$, it suffices
to take the product of the symplectic manifold $W_s$ constructed
in Proposition~\ref{2slefschetz} by a $2$--dimensional torus
$\TT^2$, and then consider a Donaldson symplectic submanifold to
reduce the dimension.

\bigskip

One can also address the problem of constructing examples of
symplectic manifolds $M_s$ which are $s$--Lefschetz and not
$(s+1)$--Lefschetz for {\em odd} integer numbers $s\geq 1$. We do
the cases $s=1$ and $s=3$. Consider the connected completely
solvable Lie group $G$ of dimension $6$ consisting of matrices of
the form
 $$
 a=\pmatrix{e^t&0&xe^t&0&0&y_1 \cr
 0&e^{-t}&0&xe^{-t}&0&y_2\cr 0&0&e^t&0&0&z_1\cr 0&0&0&e^{-t}&0&z_2
 \cr 0&0&0&0&1&t \cr 0&0&0&0&0&1 \cr},
 $$
where $t, x, y_i, z_i \in \RR$ ($i=1,2$). A global system of
coordinates $(t,x,y_1,y_2,z_1,z_2)$ for $G$ is defined by
$t(a)=t$, $x(a)=x$, $y_i(a)=y_i$, $z_i(a)=z_i$, and a standard
calculation shows that a basis for the left invariant $1$--forms
on $G$ consists of
    $$
    \{dt,\ dx,\  e^{-t}dy_1-xe^{-t}dz_1,\ e^tdy_2-xe^tdz_2,\ e^{-t}dz_1,\ e^{t}dz_2\}.
    $$
Let $\Gamma$ be a discrete subgroup of $G$ such that the quotient
space $M =\Gamma\backslash G$ is compact. (Such a subgroup exists,
see~\cite{FLS}.) Hence the forms $dt$, $dx$,
$e^{-t}dy_1-xe^{-t}dz_1$, $e^tdy_2-xe^tdz_1$, $e^{-t}dz_1$,
$e^{t}dz_2$ descend to $1$--forms $\alpha$, $\beta$, $\gamma_1$,
$\gamma_2$, $\delta_1$, $\delta_2$ on $M$ satisfying
 $$
 d\alpha=d\beta=0, \quad
 d\gamma_1=-\alpha \wedge \gamma_1 - \beta \wedge \delta_1,\quad
 d\gamma_2=\alpha \wedge \gamma_2 - \beta \wedge \delta_2,\quad
 d\delta_1=- \alpha \wedge \delta_1,\quad
 d\delta_2= \alpha \wedge \delta_2,
 $$
and such that $\{\alpha, \beta, \gamma_1, \gamma_2, \delta_1,
\delta_2 \}$ is a global basis for the $1$--forms on $M$. Using
Hattori's theorem~\cite{Hat} we compute the real cohomology of
$M$:
 \begin{eqnarray*}
 H^0(M) &=& \la 1\ra,\\
 H^1(M) &=& \la [\alpha], [\beta]\ra,\\
 H^2(M) &=& \la [\alpha \wedge \beta], [\delta_1\wedge \delta_2],
 [\gamma_1 \wedge \delta_2 + \gamma_2 \wedge \delta_1]\ra,\\
 H^3(M) &=& \la [\alpha \wedge \delta_1 \wedge \delta_2],
  [\beta\wedge \gamma_1 \wedge \gamma_2],
 [\beta\wedge (\gamma_1 \wedge \delta_2 + \gamma_2 \wedge \delta_1)],
  [\alpha\wedge (\gamma_1 \wedge \delta_2 + \gamma_2 \wedge
  \delta_1)]\ra,\\
 H^4(M) &=& \la [\alpha \wedge \beta \wedge \gamma_1 \wedge \gamma_2],
  [\alpha \wedge \beta \wedge \gamma_1 \wedge \delta_2],
   [\gamma_1 \wedge \gamma_2 \wedge \delta_1\wedge \delta_2]\ra,\\
 H^5(M) &=& \la [\alpha \wedge \gamma_1 \wedge \gamma_2 \wedge \delta_1\wedge
   \delta_2], [\beta \wedge \gamma_1 \wedge \gamma_2 \wedge \delta_1\wedge
   \delta_2]\ra, \\
 H^6(M) &=& \la [\alpha \wedge \beta \wedge \gamma_1 \wedge \gamma_2 \wedge
       \delta_1\wedge  \delta_2]\ra.
 \end{eqnarray*}
Consider the symplectic form $\omega$ on $M$ given by
$\omega=\alpha \wedge \beta + \gamma_1 \wedge \delta_2 + \gamma_2
\wedge \delta_1$. Then $[\omega] \cup [\delta_1 \wedge
\delta_2]=0$ in $H^4(M)$, which means that $M$ is not
$2$--Lefschetz. But a simple computation shows that the cup
product by $[\omega]^2$ is an isomorphism between $H^1(M)$ and
$H^5(M)$. Therefore, $(M,\omega)$ is $1$--Lefschetz, but not
$2$--Lefschetz. Moreover, $b_k^{\rm hr}(M,\omega) = b_k(M)$ for
$k\not=4$, and $b_4^{\rm hr}(M,\omega)=2<3=b_4(M)$ (compare with
Corollary \ref{6simplyconnected}). The same holds for any
symplectic form on $M$~\cite{IRTU}. Therefore, $(M,\omega)$ is
$1$--Lefschetz, but not $2$--Lefschetz.

\smallskip

Now we deal with the case $s=3$. Consider a symplectic embedding
of $(M_1,\Omega_1)=(M,\omega)$ in the complex projective space
$CP^{m_1}$, with $m_1=7$, endowed with its standard symplectic
form. We define $(M_3=\widetilde{CP}{}^{m_1},\Omega_3)$ as the
symplectic blow-up of $CP^{m_1}$ along $M_1$.

\begin{proposition}\label{2snolefschetz}
The simply connected compact symplectic manifold $(M_3,\Omega_3)$
is $3$--Lefschetz but not $4$--Lefschetz. Moreover, there is a
symplectic submanifold $W_3 \inc M_3$ of dimension~$10$ which is
$3$--Lefschetz but not $4$--Lefschetz, and every de Rham
cohomology class in $H^i(W_3)$ admits a symplectically harmonic
representative for any $i\not=6$.
\end{proposition}

\begin{proof}
Since $(M,\omega)$ is $1$--Lefschetz but not $2$--Lefschetz,
Theorem~\ref{lefschetz-blowup} and
Proposition~\ref{nolefschetz-blowup} imply that
$M_3=\widetilde{CP}{}^{7}$ is $3$--Lefschetz and not
$4$--Lefschetz. As in the proof of Proposition~\ref{2slefschetz},
an iterated Donaldson submanifold $Z_l$, $l=2$, of $M_3$ provides
an example $W_3$ in dimension $10$ which is $3$--Lefschetz and not
$4$--Lefschetz.
\end{proof}

Note also that there exists simply connected compact symplectic
manifolds of dimension~$6$ which are $1$--Lefschetz but not
$2$--Lefschetz \cite[Theorem 7.1]{Go}.

\bigskip

\noindent {\bf Acknowledgments.} This work has been partially
supported through grants MCyT (Spain) Project
BFM2001-3778-C03-02/03, UPV 00127.310-E-14813/2002, UPV
00127.310-E-15909/2004 and MTM2004-07090-C03-01.

{\small

\vspace{0.15cm}

\noindent{\sf M. Fern\'andez:} Departamento de Matem\'aticas,
Facultad de Ciencia y Tecnolog\'{\i}a, Universidad del Pa\'{\i}s
Vasco, Apartado 644, 48080 Bilbao, Spain. {\sl E-mail:}
mtpferol@lg.ehu.es

\vspace{0.15cm}

\noindent{\sf V. Mu\~noz:} Departamento de Matem\'aticas, Consejo
Superior de Investigaciones Cient{\'\i}ficas, C/ Serrano 113bis, 28006
Madrid, Spain. {\sl E-mail:} vicente.munoz@imaff.cfmac.csic.es

\vspace{0.15cm}

\noindent{\sf L. Ugarte:} Departamento de Matem\'aticas, Facultad
de Ciencias, Universidad de Zaragoza, Campus Plaza San Francisco,
50009 Zaragoza, Spain. {\sl E-mail:} ugarte@unizar.es}

\end{document}